\documentclass[11pt]{article}
\usepackage{a4wide,amsmath,amssymb,array,xspace,enumerate,bm,graphicx,color}
\begin{document}
 \newcommand{\tf}[2]{\tfrac{#1}{#2}}
\newcommand{\ST}{\texttt{ST}\xspace}
\newcommand{\tr}{\mathop{\text{trace}}}
\newcommand{\minimize}{\mathop{\text{minimize}}}
\newcommand{\maximize}{\mathop{\text{maximize}}}
\newtheorem{theorem}{Theorem}
\newtheorem{lemma}{Lemma}
\newtheorem{conjecture}{Conjecture}
\newtheorem{corollary}{Corollary}
\newcounter{remark}
\newcommand{\Real}{\mathop{\text{Re}}}
\newcommand{\remark}{\stepcounter{remark}\paragraph{Remark \theremark.}}
\newcommand{\bb}{2.4449779}
\newcommand{\bP}{\bar P}
\newcommand{\bV}{\bar V}
\newcommand{\Imag}{\mathop{\text{Im}}}
\newcommand{\diag}{\mathop{\text{diag}}}
\newcommand{\om}{\omega}
\newcommand{\AW}{{AW}\xspace}
\bibliographystyle{plain}

\stepcounter{footnote}
\author{}

\author{Richard Weber\footnote{Statistical Laboratory, Centre for
    Mathematical Sciences, Wilberforce Road, Cambridge CB2 0WB, rrw1\@@cam.ac.uk}}
\title{\LARGE \bfseries The Anderson--Weber strategy is not optimal for symmetric
  rendezvous search on {$K_4$}} \date{8 July 2009}
\newpage
\maketitle 

\begin{abstract}
We consider the symmetric rendezvous search game on a complete graph
of $n$ locations. In 1990, Anderson and Weber proposed a strategy in
which, over successive blocks of $n-1$ steps, the players
independently choose either to stay at their initial location or to
tour the other $n-1$ locations, with probabilities $p$ and $1-p$,
respectively. Their strategy has been proved optimal for $n=2$ with
$p=1/2$, and for $n=3$ with $p=1/3$. The proof for $n=3$ is very
complicated and it has been difficult to guess what might be true for
$n>3$. Anderson and Weber suspected that their strategy might not be
optimal for $n>3$, but they had no particular reason to believe this
and no one has been able to find anything better.  This paper
describes a strategy that is better than Anderson--Weber for
$n=4$. However, it is better by only a tiny fraction of a percent.
\end{abstract}


\section{The Anderson--Weber strategy}

In the symmetric rendezvous search game on $K_n$ (the completely
connected graph on $n$ vertices) two players are initially placed at
two distinct vertices (called locations). The game is played in
discrete steps and at each step each player can either stay where he
is or move to a different location. The players share no common
labelling of the locations. Our aim is to find a (randomizing)
strategy such that if both players independently follow this strategy
then they minimize the expected number of steps until they first
meet. Rendezvous search games of this type were first proposed by
Steve Alpern in 1976. They are simple to describe, and have received
considerable attention in the popular press as they model problems
that are familiar in real life. They are notoriously difficult to
analyse.

The Anderson--Weber strategy is a mixed strategy that proceeds in
blocks of $n-1$ steps.  Players begin at distinct locations, called
their home locations. In each successive block a player either stays
at his home location, with probability $p$, or makes a randomly chosen
tour of his $n-1$ non-home locations, doing this with probability
$1-p$. The motivation for the strategy comes from the wait-for-mommy
strategy that is optimal in an asymmetric version of the problem. With
probability $2p(1-p)$ the players play the wait-for-mommy strategy
over the first $n-1$ steps and so rendezvous in expected time
$(n+1)/2$.

Anderson and Weber (1990) proved that the above strategy is optimal
for the game on $K_2$, with $p=1/2$, and conjectured that it should be
optimal for $K_3$, with $p=1/3$. This was finally proved by Weber
(2006), who established a strong AW property (SAW) that \AW minimizes
$E[\min\{T,k\}]$ for all $k$. Anderson and Weber suspected that their
strategy might not be optimal for $n>3$, but they had no particular
reason to believe this and no one has been able to find any strategy
that is better. Indeed, \AW has been shown optimal amongst 2--Markov
strategies. Fan (2009) showed that \AW minimizes $P(T>2)$ and
$E[\min\{T,2\}]$. He also found that \AW is not optimal on $K_4$ if
players have the extra information that the location can be viewed as
being arranged on a circle and the players are given a common notion
of clockwise. However, the question as to whether or not \AW is
optimal has remained open for the case in which there is no such
special extra information. Fan writes, `The author believes that SAW
still holds on $K_4$, and so AW strategy is still optimal'. We were
inclined to agree, but now find that \AW can be bettered. For more
background to the problem see Weber (2006).

Let us begin by reprising the \AW strategy for the symmetric
rendezvous game on $K_4$. We assume that there is no special knowledge
(such as a common notion of clockwise on a circle).  The \AW strategy
is a 3--Markov strategy that repeats in blocks of 3 steps. In each
successive block of 3 steps, each player, independently, remains at
his home location with probability $p$, or does a random chosen tour
of his 3 non-home locations, with probability $1-p$. This leads to
rendezvous in an expected number steps $ET$, where
\begin{align*}
ET &= p^2\times (3+ET)  + 2p(1-p)\times 2 
+ (1-p)^2\times \left(\tfrac{1}{2}(16/9)+\tfrac{1}{2}(3+ET)\right)\\
&=\frac{43-14 p + 25 p^2}{9 \left(1+2 p-3 p^2\right)}\,.
\end{align*}
This is explained as follows.
\begin{enumerate}
\item If both stay home they do not meet.
\item If one stays home, while the other tours, then they meet in
  expected time 2.
\item If both tour, then they meet with probability $1/2$, and
  conditional on meeting they meet in expected time $16/9$.
\end{enumerate}
One easily finds that the minimum of $ET$ is achieved by taking
\[
p=\frac{1}{4} \left(3 \sqrt{681}-77\right)\approx 0.321983
\]
and then
\[
ET=\frac{1}{12} \left(15+\sqrt{681}\right)\approx 3.42466\,.
\]\medskip

\section{A strategy better than Anderson--Weber on $K_4$}

We now explain how the \AW strategy can be bettered.  Suppose player I
has location 1 as his home, and player II has location 2 as home.  We
might imagine that each player labels his non-home locations as
$a,b,c$, and so a tour of his non-home locations is one of six
possible tours: $abc$, $acb$, $bac$, $bca$, $cab$, $cba$. In the case that
player I has $(a,b,c)=(2,3,4)$ and player II has $(a,b,c)=(1,3,4)$ we
can compute the matrix
\[
B=\left(
\begin{array}{llllll}
 2 & \text{X} & 3 & \text{X} & \text{X} & 2 \\
 \text{X} & 2 & \text{X} & 2 & 3 & \text{X} \\
 3 & \text{X} & 1 & 1 & \text{X} & \text{X} \\
 \text{X} & 2 & 1 & 1 & \text{X} & \text{X} \\
 \text{X} & 3 & \text{X} & \text{X} & 1 & 1 \\
 2 & \text{X} & \text{X} & \text{X} & 1 & 1
\end{array}
\right)
\]
where we have ordered the rows and columns to correspond to $abc$,
$acb$, $bac$, $bca$, $cab$, $cba$. A number entry indicates the step
at which players meet when they meet, and X indicates that they do not
meet. There are 36 such matrices, over which we must average, for each
possible pair of assignments by players I and II, of $(2,3,4)$ and
$(1,3,4)$, respectively, to $(a,b,c)$.

Let us begin by noting that if a player stays home for three steps and
meeting does not occur, then the other player must also have been
staying home. Similarly, if a player tours for three steps and meeting
does not occur, then the other player must also have been touring (and
their tours not meeting). Thus after any $3k$ steps (a multiple of 3)
each player knows exactly how many times both have been touring.

Whenever a player makes a tour in the \AW strategy he chooses his tour
at random (independently of previous tours). We show how to improve
\AW introducing some dependence between tours. Let us adopt a notation
in which the first tour a player makes is labelled $A$. The second
distinct tour a player makes is labelled $B$, and so on. So, for
example, $AAB$ means that on his first three tours, a player (i) first
makes a random tour, (ii) second makes the same tour as his first tour,
(iii) and third makes a tour chosen randomly from amongst the 5 tours
he has not yet tried.

Let us consider first a modified problem in which at each so-called
`t--step' each player makes a tour of his non-home locations. In this
modified problem no player stays home for a t--step.  We wish to
minimize the expected number of t--steps until the players meet. At
the first t--step both players do $A$ and the probability of meeting
is $1/2$.  If a 1--Markov strategy is employed, so successive t--steps
are chosen at random, then the expected number of t--steps until meeting
occurs is 2.

Over the first two t--steps, the players can do either $AA$ or $AB$. The
matrix for not meeting is
\[
P_2=\left(
\begin{array}{cc}
 \frac{1}{2} & \frac{1}{5} \\[3pt]
 \frac{1}{5} & \frac{13}{50}
\end{array}
\right)
\]
One can check that $P_2\succ 0$ (i.e., $P_2$ is positive definite).
thus for a 2--Markov strategy we would be solving
\[
ET = x^\top(J+P_2\, ET)x
\]
where $J$ is a $2\times 2$ matrix filled with 1s.  This has a minimum
value of $ET=2$, when we take $x^\top=(1/6,5/6)$. This means that,
restricting to 2--Markov strategies, tours should be chosen at
random. 

Similarly, over the first three t--steps, the players can do $AAA$,
$AAB$, $ABA$, $ABB$, $ABC$. The matrix for not meeting is
\[
P_3=\left(
\begin{array}{ccccc}
 \frac{1}{2} & \frac{1}{5} & \frac{1}{5} & \frac{1}{5} & \frac{1}{20} \\[3pt]
 \frac{1}{5} & \frac{13}{50} & \frac{2}{25} & \frac{2}{25} & \frac{11}{100} \\[3pt]
 \frac{1}{5} & \frac{2}{25} & \frac{13}{50} & \frac{2}{25} & \frac{11}{100} \\[3pt]
 \frac{1}{5} & \frac{2}{25} & \frac{2}{25} & \frac{13}{50} & \frac{11}{100} \\[3pt]
 \frac{1}{20} & \frac{11}{100} & \frac{11}{100} & \frac{11}{100} & \frac{7}{50}
\end{array}
\right)\,.
\]
Again, $P_3\succeq 0$, and $(x_{AAA}, x_{AAB}, x_{ABA}, x_{ABB},
x_{ABC})=(1/36,5/36,5/36,5/36,20/36)$ is optimal in the sense of
minimizing the solution of
\[
ET = x^\top(J + P_2 + P_3\, ET)x\,,
\]
where $J$ is now $5\times 5$ and $P_2$ is expanded to the appropriate
$5\times 5$ matrix. Thus
amongst 3--Markov strategies, tours should also be chosen at random. 

However, over four t--steps things turn out differently. There are now
15 possible strategies: $AAAA$, $AAAB$, $AABA$, $AABB$, $AABC$,
$ABAA$, $ABAB$, $ABAC$, $ABBA$, $ABBB$, $ABBC$, $ABCA$, $ABCB$,
$ABCC$, $ABCD$. The matrix for not meeting can be computed to be

\[
P_4=\left(
\begin{array}{ccccccccccccccc}
 \frac{1}{2} & \frac{1}{5} & \frac{1}{5} & \frac{1}{5} & \frac{1}{20} & \frac{1}{5} &
   \frac{1}{5} & \frac{1}{20} & \frac{1}{5} & \frac{1}{5} & \frac{1}{20} & \frac{1}{20} &
   \frac{1}{20} & \frac{1}{20} & 0 \\[3pt]
 \frac{1}{5} & \frac{13}{50} & \frac{2}{25} & \frac{2}{25} & \frac{11}{100} & \frac{2}{25} &
   \frac{2}{25} & \frac{11}{100} & \frac{2}{25} & \frac{2}{25} & \frac{11}{100} & \frac{1}{50}
   & \frac{1}{50} & \frac{1}{50} & \frac{3}{100} \\[3pt]
 \frac{1}{5} & \frac{2}{25} & \frac{13}{50} & \frac{2}{25} & \frac{11}{100} & \frac{2}{25} &
   \frac{2}{25} & \frac{1}{50} & \frac{2}{25} & \frac{2}{25} & \frac{1}{50} & \frac{11}{100} &
   \frac{11}{100} & \frac{1}{50} & \frac{3}{100} \\[3pt]
 \frac{1}{5} & \frac{2}{25} & \frac{2}{25} & \frac{13}{50} & \frac{11}{100} & \frac{2}{25} &
   \frac{2}{75} & \frac{1}{30} & \frac{2}{75} & \frac{2}{25} & \frac{1}{30} & \frac{1}{30} &
   \frac{1}{30} & \frac{11}{100} & \frac{23}{450} \\[3pt]
 \frac{1}{20} & \frac{11}{100} & \frac{11}{100} & \frac{11}{100} & \frac{7}{50} & \frac{1}{50}
   & \frac{1}{30} & \frac{7}{150} & \frac{1}{30} & \frac{1}{50} & \frac{7}{150} & \frac{7}{150}
   & \frac{7}{150} & \frac{1}{20} & \frac{14}{225} \\[3pt]
 \frac{1}{5} & \frac{2}{25} & \frac{2}{25} & \frac{2}{25} & \frac{1}{50} & \frac{13}{50} &
   \frac{2}{25} & \frac{11}{100} & \frac{2}{25} & \frac{2}{25} & \frac{1}{50} & \frac{11}{100}
   & \frac{1}{50} & \frac{11}{100} & \frac{3}{100} \\[3pt]
 \frac{1}{5} & \frac{2}{25} & \frac{2}{25} & \frac{2}{75} & \frac{1}{30} & \frac{2}{25} &
   \frac{13}{50} & \frac{11}{100} & \frac{2}{75} & \frac{2}{25} & \frac{1}{30} & \frac{1}{30} &
   \frac{11}{100} & \frac{1}{30} & \frac{23}{450} \\[3pt]
 \frac{1}{20} & \frac{11}{100} & \frac{1}{50} & \frac{1}{30} & \frac{7}{150} & \frac{11}{100} &
   \frac{11}{100} & \frac{7}{50} & \frac{1}{30} & \frac{1}{50} & \frac{7}{150} & \frac{7}{150}
   & \frac{1}{20} & \frac{7}{150} & \frac{14}{225} \\[3pt]
 \frac{1}{5} & \frac{2}{25} & \frac{2}{25} & \frac{2}{75} & \frac{1}{30} & \frac{2}{25} &
   \frac{2}{75} & \frac{1}{30} & \frac{13}{50} & \frac{2}{25} & \frac{11}{100} & \frac{11}{100}
   & \frac{1}{30} & \frac{1}{30} & \frac{23}{450} \\[3pt]
 \frac{1}{5} & \frac{2}{25} & \frac{2}{25} & \frac{2}{25} & \frac{1}{50} & \frac{2}{25} &
   \frac{2}{25} & \frac{1}{50} & \frac{2}{25} & \frac{13}{50} & \frac{11}{100} & \frac{1}{50} &
   \frac{11}{100} & \frac{11}{100} & \frac{3}{100} \\[3pt]
 \frac{1}{20} & \frac{11}{100} & \frac{1}{50} & \frac{1}{30} & \frac{7}{150} & \frac{1}{50} &
   \frac{1}{30} & \frac{7}{150} & \frac{11}{100} & \frac{11}{100} & \frac{7}{50} & \frac{1}{20}
   & \frac{7}{150} & \frac{7}{150} & \frac{14}{225} \\[3pt]
 \frac{1}{20} & \frac{1}{50} & \frac{11}{100} & \frac{1}{30} & \frac{7}{150} & \frac{11}{100} &
   \frac{1}{30} & \frac{7}{150} & \frac{11}{100} & \frac{1}{50} & \frac{1}{20} & \frac{7}{50} &
   \frac{7}{150} & \frac{7}{150} & \frac{14}{225} \\[3pt]
 \frac{1}{20} & \frac{1}{50} & \frac{11}{100} & \frac{1}{30} & \frac{7}{150} & \frac{1}{50} &
   \frac{11}{100} & \frac{1}{20} & \frac{1}{30} & \frac{11}{100} & \frac{7}{150} &
   \frac{7}{150} & \frac{7}{50} & \frac{7}{150} & \frac{14}{225} \\[3pt]
 \frac{1}{20} & \frac{1}{50} & \frac{1}{50} & \frac{11}{100} & \frac{1}{20} & \frac{11}{100} &
   \frac{1}{30} & \frac{7}{150} & \frac{1}{30} & \frac{11}{100} & \frac{7}{150} & \frac{7}{150}
   & \frac{7}{150} & \frac{7}{50} & \frac{14}{225} \\[3pt]
 0 & \frac{3}{100} & \frac{3}{100} & \frac{23}{450} & \frac{14}{225} & \frac{3}{100} &
   \frac{23}{450} & \frac{14}{225} & \frac{23}{450} & \frac{3}{100} & \frac{14}{225} &
   \frac{14}{225} & \frac{14}{225} & \frac{14}{225} & \frac{7}{90}
\end{array}
\right)
\]
It now turns out that $P_4$ has a negative eigenvalue. The \AW
strategy would be to choose tours at random, which gives
\begin{align*}
x^\top &=
\Bigl(p_{AAAA}, p_{AAAB}, p_{AABA}, p_{AABB}, p_{AABC},
p_{ABAA}, p_{ABAB}, p_{ABAC}, p_{ABBA}, p_{ABBB}, \\
&p_{ABBC}, p_{ABCA}, p_{ABCB},
p_{ABCC}, p_{ABCD}\Bigr)\\
&= \frac{1}{6^3}(1, 5, 5, 5, 20, 5, 5, 20, 5, 5, 20, 20, 20, 20, 60)\,.
\end{align*}
Solving $ET=x^\top(J+P_2+P_3+P_4\, ET)x$, we find $ET=2$, as we
expect. However consider
\[
y^\top=\left(0,\frac{1}{12},\frac{1}{12},0,0,\frac{1}{12},0,0,0,\frac{1}{12},0,0,0,0,\frac{2}{3}\right)\,.
\]
Solving $ET=y^\top(J+P_2+P_3+P_4\, ET)y$ gives
$ET=2-\frac{23}{16200}=1.99858$. Thus, rendezvous occurs in a smaller
expected number of t--steps than it does under \AW. This happens when players
use a mixed 4--Markov strategy of doing $AAAB$, $AABA$, $ABAA$,
$ABBB$, each with probability $1/12$, and $ABCD$ with probability
$2/3$. This corresponds to choosing tours for the first two t--steps
at random, but then making the choice of tours at the 3rd and 4th
t--step depend on the tours taken at the 1st and 2nd t--step. The
choice of $y$ is not unique. It has been choose to be simple,
containing many 0s, and it was found by using the fact that the
eigenvector of $P_4$ having a negative eigenvalue is of a pattern
$(\alpha,\beta,\beta,\gamma,\delta,\beta,\gamma,\delta,\gamma,\beta,\delta
,\delta,\delta,\delta,\epsilon)$ for some irrational $\alpha$,
$\beta$, $\gamma$, $\delta$, $\epsilon$.  \medskip

The above makes it very plausible that we can find a strategy that is
better than \AW on $K_4$. We now need to do some careful calculations.
We consider a 12--Markov strategy consisting of 4 t--steps. In each
t-step a player remains home with probability $p$, and tours with
probability $1-p$. When he makes tours, he does so in an manner that
achieves the distribution previously described. That is, any 1st and
2nd tours are made at random, but 3rd and 4th tours are made so that
these are consistent with the distribution over 4 tours being $AAAB$,
$AABA$, $ABAA$, $ABBB$, each with probability $1/12$, and $ABCD$ with
probability $2/3$. If at the end of 12 steps the players have not met
then the strategy restarts, forgetting about the number of previous
t-steps on which players made non-meeting tours.

We found it easiest to calculate the expected meeting time by
attaching a probability to each possible 12--step paths that the
strategy might take. There are 1585 possible paths which have nonzero
probability. We computed the step at which players meet, or event
that they do not meet, for each of the $1585\times 1585$
possibilities, and averaged these using the appropriate
probabilities. The calculations are intricate, but can be checked in
various ways to provide confidence that no mistake has been made.
It turns out that the expected meeting time is

\scriptsize
\[
\text{\normalsize $ET=$}{\frac{-227773 p^8+582884 p^7-1329319 p^6+1737938 p^5-1941235 p^4+1420688 p^3-998569 p^2+389834
   p-217648}{3 \left(82001 p^8-218608 p^7+327728 p^6-315256 p^5+215870 p^4-104656 p^3+36128
   p^2-8008 p-15199\right)}}\,.
\]

\normalsize Taking $p=\frac{1}{4} \left(3 \sqrt{681}-77\right)$, which
is the optimal value for the \AW strategy, we find that the new
strategy produces an expected meeting time that is less than that of
\AW by
\[
\frac{243 \left(75041961207+4700853101
    \sqrt{681}\right)}{327540887401488016}\approx 0.000146683\,.
\]
The tiny improvement is due to the fact that when both players do four
t-tours (which happens with probability $(1-p)^4$), the new strategy
gives a greater probability that the players meet than does \AW.  It
would be possible to make the new strategy even better, by choosing
$p$ slightly differently, or indeed making it depend on the number of
tours that have been taken so far over which players have not met. We
could also do better by not restarting after 12 steps. However, our
aim is not to try to find the best strategy for $K_4$, which still seems
very difficult, but simply to show that \AW is not optimal. This
we have now done.

\end{document}